\documentclass[reqno,12pt]{amsart}

\usepackage{epsf}
\usepackage{graphics}
\usepackage{amssymb}
\usepackage{amsmath}

\date{}

\theoremstyle{plain}
\newtheorem{theorem}{Theorem}
\newtheorem{corollary}{Corollary}

\theoremstyle{definition}

\theoremstyle{remark}

\newtheorem*{remarks}{Remarks}

\def\N{{\mathbb N}}

\def\R{{\mathbb R}}

\title{M\"obius Bands with a Quasipositive Fibred Hole} 

\author{Sebastian Baader}

\begin{document}

\begin{abstract} We prove that every knot in the 3-space bounds an embedded punctured M\"obius band whose other boundary component is a quasipositive fibred knot.
\end{abstract}

\maketitle

\section{Introduction} 

Quasipositive braids are products of conjugates of a positive standard
generater of the braid group; they generalize the class of positive braids. A
quasipositive knot is a knot which is isotopic to the closure of a
quasipositive braid. Quasipositive knots emerged from the theory of complex
plane curves (\cite{Ru2}). From certain viewpoints, quasipositive knots are
`dense' in the set of knots. For instance, every Alexander polynomial can be
realized by a quasipositive knot (\cite{Ru1}). A similar statement is true for
Vassiliev invariants (\cite{Ba}). In this paper, we prove another density
result for quasipositive knots, involving unoriented cobordisms. A cobordism between two disjoint links $L_1$, $L_2 \subset \R^3$ is a smooth embedded surface $S \subset \R^3$ whose boundary is precisely $L_1 \cup L_2$. We distinguish between oriented and unoriented cobordisms, depending on whether the surface and its boundary components are oriented or not. In this paper, we consider particular unoriented cobordisms, namely punctured M\"obius bands embedded in $\R^3$.

\begin{theorem}
For every knot $K \subset \R^3$ there exists a quasipositive fibred knot $Q \subset \R^3$ and an embedded punctured M\"obius band $M \subset \R^3$ whose boundary components are precisely $K$ and $Q$.
\end{theorem}

\section{From Braids to Quasipositive Braids}

Let $\beta \in B_n$ be an element of the braid group on $n$ strings. We may
assume that the first letter of $\beta$ is a negative standard generator of
$B_n$, since we can always insert a canceling pair $\sigma_1^{-1} \sigma_1$
in front of $\beta$. Therefore, there exist a natural number $k \in \N$,
natural numbers $i_1, i_2, \ldots, i_k \in \{1, 2, \ldots, n-1\}$, and positive
(possibly empty) words $w_1, w_2, \ldots, w_k \in B_n$, such that
\begin{equation}
\beta=\prod^{k}_{l=1} \sigma_{i_l}^{-1} w_l.
\label{beta}
\end{equation}
Let $K$ be an arbitrary knot. $K$ is isotopic to the closure of a braid $\beta
\in B_n$, written as in (\ref{beta}). The closure of $\beta$ is depicted in
figure~1.

\begin{figure}[ht]
\scalebox{0.8}{\raisebox{-0pt}{$\vcenter{\hbox{\epsffile{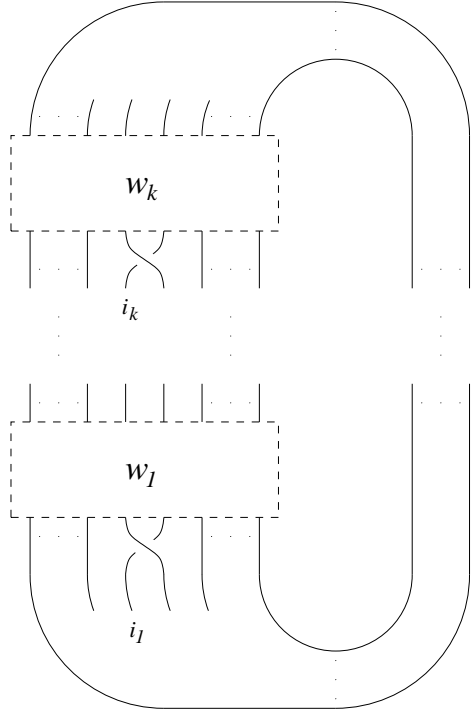}}}$}}
\caption{Closure of the braid $\beta$}
\end{figure}

Now, we construct a new braid $\alpha \in B_{n+2}$, as follows:
\begin{equation}
\aligned
\alpha= \sigma_{n+1} \sigma_{n} \prod^{k}_{l=1} 
&(\sigma_{n} \sigma_{n-1} \cdots \sigma_{i_l} 
\sigma_{n+1} \sigma_{n} \cdots \sigma_{i_l+1}
\sigma_{i_l+2}^{-1} \\
&\quad 
\sigma_{i_l+1} \sigma_{i_l+2} \cdots \sigma_{n+1} 
\sigma_{i_l} \sigma_{i_l+1} \cdots \sigma_{n}
w_l).
\endaligned
\label{alpha}
\end{equation}

$\alpha$ is clearly a quasipositive braid, since all its negative letters
appear in conjugating pairs, together with their positive counterparts:
$$\sigma_{i_l+2}^{-1} \sigma_{i_l+1} \sigma_{i_l+2}.$$ 
Let $Q$ be the closure of $\alpha$, as depicted in figure~2. A quick
comparison with figure~1 shows that $Q$ is also a knot. In particular, $Q$ is
a quasipositive knot. 

\begin{figure}[ht]
\scalebox{0.8}{\raisebox{-0pt}{$\vcenter{\hbox{\epsffile{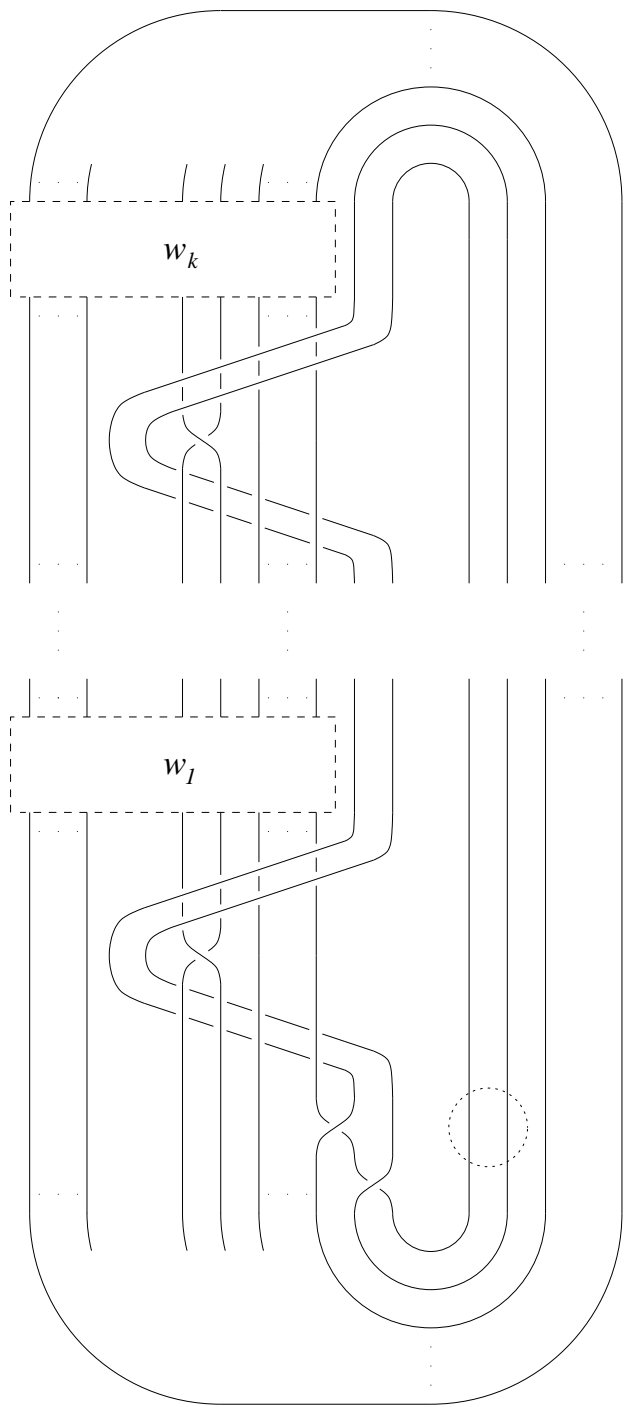}}}$}}
\caption{Closure of the braid $\alpha$}
\end{figure}

Next, we construct an embedded punctured M\"obius band whose boundary components are $K$ and $Q$: choose any embedded annulus $A \subset \R^3$, one of whose boundary components is the closure of $\alpha$, i.e. $K$. Glue a M\"obius band to $A$, woven as the additional band of $\beta$, as shown in figure~3. The two boundary components of the resulting punctured M\"obius band are isotopic to $K$ and $Q$.

\begin{figure}[ht]
\scalebox{1}{\raisebox{-0pt}{$\vcenter{\hbox{\epsffile{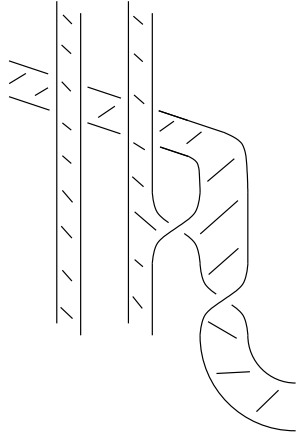}}}$}}
\caption{Gluing of M\"obius band}
\end{figure}

At last, we construct a natural fibre surface $F$ for $Q$. $F$ is a union of $(n+2)$ discs, one for each braid strand of 
$\alpha$, and several bands with a positive half twist. Each positive crossing corresponds to such a twisted band, except in the neighbourhood of negative crossings, where we choose twisted bands as shown in figure~4.

\begin{figure}[ht]
\scalebox{1}{\raisebox{-0pt}{$\vcenter{\hbox{\epsffile{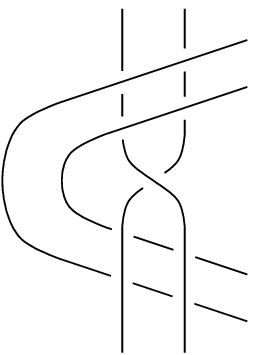}}}$}} \qquad \qquad
\scalebox{1}{\raisebox{-0pt}{$\vcenter{\hbox{\epsffile{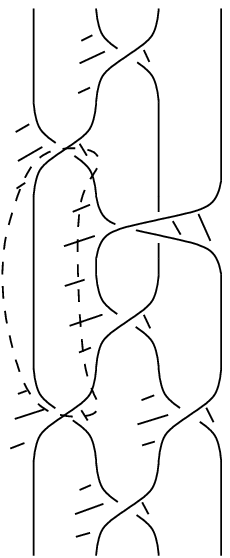}}}$}}
\caption{Seifert surface $F$ for $Q$}
\end{figure}

The dashed line indicates the core curve of a Hopf band in $F$, which we may deplumb without affecting the attribute of fibredness or non-fibredness. Thus we may cut the upper left twisted band and obtain the surface $F'$ shown on the left hand side of figure~5. 
The surface $F'$ in turn is isotopic to the canonical Seifert surface of a positive braid, as shown on the right hand side of figure~5 (slide the long band along the upper band). This concludes the proof of theorem~1, since canonical Seifert surfaces of positive braids are fibre surfaces.

\begin{figure}[ht]
\scalebox{1}{\raisebox{-0pt}{$\vcenter{\hbox{\epsffile{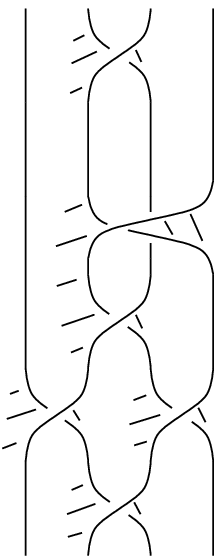}}}$}} \qquad \qquad
\scalebox{1}{\raisebox{-0pt}{$\vcenter{\hbox{\epsffile{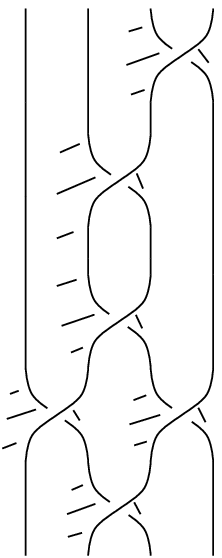}}}$}}
\caption{Seifert surface $F'$ after deplumbing}
\end{figure}

\begin{remarks} \quad \\
\begin{enumerate}
\item The knot $Q$ is even strongly quasipositive, in the sense of L.~Ru\-dolph (\cite{Ru3}). \\

\item It is absolutely essential that the surface $M$ be orientable. Indeed, there exist strong restrictions on the Euler number of oriented cobordisms between quasinegative and quasipositive knots (e.g. the Rasmussen invariant (\cite{Ra})). \\

\item The application of one saddle point move inside the dashed circle at the bottom right of figure~2 transforms 
$Q$ into $K$. Here a saddle point move is a local move that acts on link diagrams, as shown in figure~6. In \cite{HNT}, Hoste, Nakanishi and Taniyama proved that a simple crossing change can be expressed as a sequence of two saddle point moves. Therefore, it is an unknotting operation.\\
\end{enumerate}
\end{remarks}

\begin{figure}[ht]
\scalebox{1}{\raisebox{-0pt}{$\vcenter{\hbox{\epsffile{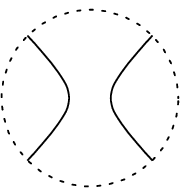}}}$}} \qquad
$\longleftrightarrow$ \qquad 
\scalebox{1}{\raisebox{-0pt}{$\vcenter{\hbox{\epsffile{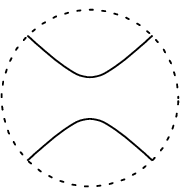}}}$}}
\caption{Saddle point move}
\end{figure}

\begin{corollary} Every knot can be transformed into a quasipositive fibred knot by the application of one saddle point move.
\end{corollary}

\bigskip
\noindent
Department of Mathematics,
ETH Z\"urich, 
Switzerland

\bigskip
\noindent
\emph{sebastian.baader@math.ethz.ch}

\end{document}